\documentclass[12pt]{article}

\usepackage{amssymb, amsthm, amsmath, color, fontenc}
\usepackage{epsfig}

\newtheorem{thm}[subsection]{Theorem}
\newtheorem{cor}[subsection]{Corollary}
\newtheorem{lem}[subsection]{Lemma}

\title{Examples of relative deformation spaces that are not locally connected}
\author{Aaron Magid \footnote{The author was partially supported by the NSF RTG grant \#0602191.}}
\date{13 April 2008}

\begin{document}
\maketitle

\begin{abstract}  
We provide an infinite family of pared manifolds whose relative deformation spaces of hyperbolic structures on these manifolds are not locally connected.  This is a natural extension of the recent result of Bromberg that shows the space of Kleinian punctured torus groups is not locally connected.   
\end{abstract}

%%%%
%%%% SECTION about introduction and statement of results
%%%%
\section{Introduction}  

One of the most recent developments in bumponomics, the study of the pathological topological properties of deformation spaces of hyperbolic 3-manifolds, was Ken Bromberg's result that the space of Kleinian punctured torus groups is not locally connected \cite{B}.  We extend this to show that there are infinitely many pared manifolds whose relative deformation spaces are not locally connected.  

Given a compact, orientable 3-manifold $M$ with incompressible boundary, the deformation space, $AH(M)$, is the set of equivalence classes of marked hyperbolic 3-manifolds homotopy equivalent to $M$.  We equip this space with the algebraic topology.  The components of the interior of this space are in one-to-one correspondence with the marked homeomorphism types of 3-manifolds homotopy equivalent to $M$.  Points in the interior correspond to geometrically finite and minimally parabolic hyperbolic manifolds. Using quasiconformal deformations, one can parameterize each of these components by the Teichm\"uller space of the conformal boundary.  This parameterization combines work of Ahlfors, Bers, Kra, Marden, Maskit, Sullivan, and Thurston.  All of this work extends to the pared manifold setting where one considers the relative deformation space, $AH(M,P)$, for some disjoint collection, $P$, of tori and annuli in $\partial M$. This is the space of hyperbolic 3-manifolds homotopy equivalent to $M$ with cusps associated to each component of $P$.  More precise statements will be given in the following section.  See Chapter 7 of \cite{CM} for a more complete discussion. 

Points in $AH(M)$ that are not in the interior are less well-behaved.  The term bumponomics was coined after Anderson, Canary, and McCullough showed that components of the interior of $AH(M)$ can bump, (\textit{i.e.}, there exist points in the closure of multiple components) \cite{AC}, \cite{ACM}.  Under certain conditions, a component $B$ of the interior of $AH(M)$ can also self-bump.  This means there is a point $\rho \in \overline{B}$ such that for every sufficiently small neighborhood $U$ of $\rho$, $U \cap B$ is disconnected.  McMullen was the first to observe that when $S$ is a surface, the space of quasi-Fuchsian representations of $\pi_1(S)$, equivalently the interior of $AH(S \times I)$, self-bumps \cite{McMullen}.  
Bromberg and Holt generalized this result to show that if $M$ contains a primitive, essential, boundary incompressible annulus whose core curve is not homotopic into a torus component of $\partial M$ then the interior of $AH(M)$ self-bumps \cite{BH}.  

 If we let $\hat{T}$ denote the punctured torus, Bromberg's result that $AH(\hat{T} \times I, \partial \hat{T} \times I)$ is not locally connected shows that self-bumping may be considerably more complicated than once thought.   He conjectured that for any compact surface $S$, $AH(S \times I, \partial S \times I)$ is not locally connected.  Although not resolving this conjecture, we provide an infinite collection of pared manifolds $(M,P)$ for which $AH(M,P)$ is not locally connected. In particular, we prove the following.

\begin{thm} \label{main theorem}
Let $M$ be a hyperbolizable 3-manifold containing a primitive essential annulus $A$, and suppose $(\hat{T} \times I, \partial \hat{T} \times I)$ is pared homeomorphic to $(M', A)$, where $M'$ is the closure of one of the components of $M-A$. If $P \subset \partial M$ is a paring locus that contains exactly one of the components of $\partial A$ and is otherwise disjoint from $M'$, then $AH(M,P)$ is not locally connected. 
\end{thm}

We remark that the hypotheses of the main theorem are not vacuous and the simplest examples of manifolds $(M,P)$ for which $AH(M,P)$ is not locally connected are described in the following corollary.  These relative deformation spaces arise naturally as subsets of the boundaries of spaces of quasi-Fuchsian representations.

\begin{cor}
Let $S$ be a closed surface of genus $g \geq 2$. Let $M = S \times I$ and $P$ a single annulus on $S \times \{1\}$ such that $P$ separates $S \times \{1\}$ into a punctured torus and a once-punctured genus $(g-1)$ surface.  Then $AH(S \times I, P)$ is not locally connected.
\end{cor}

We now outline the rest of the paper.  After reviewing some background and necessary notation, we review Bromberg's results in Section \ref{bromberg section}.  We then prove a topological lemma that detects a lack of local connectivity.  The main construction appears in Section \ref{main proof}.  Any representation in $AH(M,P)$ restricts to a representation in $AH(\hat{T} \times I, \partial \hat{T} \times I)$, so we obtain a continuous map $\Pi: AH(M,P) \to AH(\hat{T} \times I, \partial \hat{T} \times I)$.  
In order to show $AH(M,P)$ is not locally connected, we consider a sequence $\{ \rho_n \} \subset AH(\hat{T} \times I, \partial \hat{T} \times I)$ 
that converges to a point $\rho_\infty$ where 
$AH(\hat{T} \times I, \partial \hat{T} \times I)$ is not locally connected.  The construction of this sequence utilizes Bromberg's description of $AH(\hat{T} \times I, \partial \hat{T} \times I)$.  We then define a continuous section on this sequence.  In particular, if $M_0$ denotes the component of $M-A$ that is not homeomorphic to $\hat{T} \times I$, then we find a fixed representation $\eta_0 \in AH(M_0, A)$ and use Klein-Maskit combination  to construct
 representations $ \eta_n \in AH(M,P)$ from $\eta_0$ and $\rho_n$ such that $\Pi(\eta_n) = \rho_n$.  We show $\eta_n$ converges to a representation $\eta_\infty \in AH(M,P)$ such that $\Pi(\eta_\infty) = \rho_\infty$. Our topological lemma from Section \ref{topological section} then implies $AH(M,P)$ is not locally conntected at $\eta_\infty$.

Using a similar construction, we also show the following. 

\begin{thm}  \label{second theorem}
Suppose $M$ is obtained from a boundary connected sum of $\hat{T} \times I$ and a hyperbolizable manifold $(M',P')$ by identifying a disk in $\hat{T} \times \{ 0 \}$ to a disk in $\partial M' - P'$.  Let $P = P' \cup (\partial \hat{T} \times I)$. Then $AH(M,P)$ is not locally connected. 
\end{thm}

 The proof is similar to the proof of Theorem \ref{main theorem} and we also have an infinite family of examples of manifolds satisfying the hypotheses of Theorem \ref{second theorem}.   Since $\hat{T} \times I$ is homeomorphic to a genus 2 handlebody, we can form any genus $g \geq 2$ handlebody, $H$, by taking a boundary connected sum of $\hat{T} \times I$ and a genus $g-2$ handlebody $M'$. Then clearly $(H, \partial \hat{T} \times I)$ satisfies the hypotheses giving us the following corollary. 
  
\begin{cor}
Let $H$ be a genus $g\geq 2$ handlebody.  There exists an annulus $P \subset \partial H$ such that $AH(H,P)$ is not locally connected.
\end{cor}

We point out this straightforward corollary as $AH(H,P)$ arises as a subset of the boundary of the space of genus $g$ Schottky groups.  

I would like to thank Richard Canary and Ken Bromberg for helpful conversations and instructions.  Without them, this paper would not exist.

%%%%
%%%% SECTION aboutbackground
%%%%
\section{Background and Definitions}   

\subsection{Deformation Spaces}
Let $(M,P)$ be a hyperbolizable pared 3-manifold. Recall this means that $M$ is a compact, oriented, 3-manifold that is not a 3-ball, and $P \subset \partial M$ is a disjoint collection of incompressible annuli and tori, containing all tori in $\partial M$, and satisfying the property that every $\pi_1$-injective map $(S^1 \times I, S^1 \times \partial I) \to (M,P)$ is homotopic (as a map of pairs) into $P$.

We define the relative character variety 
\[
R(M,P) = Hom_P (\pi_1(M), PSL(2, \mathbb{C}))//PSL(2, \mathbb{C})
\] to be the set of conjugacy classes of representations $\rho: \pi_1(M) \to PSL(2, \mathbb{C})$ such that $\rho(g)$ is parabolic or the identity whenever $g \in \pi_1(P)$.  Let $AH(M,P)$ denote the subset of $R(M,P)$ consisting of discrete and faithful representations.  Thus $AH(M,P)$ inherits a topology from the character variety known as the algebraic topology.  Results of Chuckrow \cite{Chuckrow} and J\o rgensen \cite{Jorgensen} show that $AH(M,P)$ is a closed subset of $R(M,P)$ with respect to this topology.  See Chapter 5 of $\cite{CM}$ for more details. 

Results of Marden \cite{Marden open} and Sullivan \cite{Sullivan} show that the interior of $AH(M,P)$ consists of precisely the geometrically finite representations for which $\rho(g)$ is parabolic if and only if $g \in \pi_1(P)$. We call these minimally parabolic and denote the set of minimally parabolic representations by $MP(M,P)$.

If $G_n$ is a sequence of Kleinian groups, then we say $G_n$ converges geometrically to a Kleinian group $G$ if $G_n$ converges to $G$ with respect to the Chabauty topology on closed subsets of $PSL(2, \mathbb{C})$.  Suppose $G_n$ is a sequence of Kleinian groups converging geometrically to a Kleinian group $G$, and $\rho_n \to \rho$ is a convergent sequence in $Hom_P (\pi_1(M), PSL(2, \mathbb{C}))$ with $\rho_n(\pi_1(M)) = G_n$.  Then $\rho(\pi_1(M)) \subset G$; however, as is often the case in this paper, the containment can be strict.  If $N_n = \mathbb{H}^3/G_n$ and $N = \mathbb{H}^3/G$ are the associated hyperbolic 3-manifolds then we say $N_n$ converges geometrically to $N$. (See \cite{JM} or \cite{Marden} for more on geometric convergence.)

 %%%%%%  drilling and filling background section
 \subsection{Hyperbolic Drilling and Filling} \label{filling notation}
 Given a hyperbolic 3-manifold $N$ and simple closed geodesic $\gamma \subset N$, let $W$ be an open tubular neighborhood of $\gamma$.  Define $\hat{N} = N - W$ to be the manifold obtained by drilling out $\gamma$.  If $\phi: \hat{N} \to N$ denotes the inclusion map, then we give $\hat{N}$ the unique complete hyperbolic metric so that $\phi$ extends to a conformal map between the conformal boundaries of $N$ and $\hat{N}$.  The existence of such a metric is given in \cite{Kojima} and we say that $\hat{N}$ is the $\gamma$-drilling of $N$. 
 
To fill a hyperbolic manifold with a rank-$2$ cusp, we need a bit more notation.  Let $\hat{M}$ be a compact 3-manifold with a single torus boundary component $T$, and choose generators $m$ and $l$ for the fundamental group of this boundary component so $\pi_1(T) = \left< m \right> \oplus \left< l \right>$.  Let $M(p,q)$ denote the $(p,q)$-Dehn filling of $\hat{M}$.  By this we mean the result of gluing a solid torus $V$ to $\hat{M}$ by a homeomorphism identifying $\partial V$ with $T$ that takes the meridian of $\partial V$ to a curve in the homotopy class of $pm + ql \in \pi_1(T)$.  

Let $\hat{N}$ denote the interior of $\hat{M}$ with a fixed, complete, geometrically finite hyperbolic structure.  Then, if it exists, we let $N(p,q)$ denote the interior of $M(p,q)$ together with a complete hyperbolic structure such that the inclusion map $\hat{N} \to N(p,q)$ extends to a conformal map between the conformal boundary components.  When such a hyperbolic metric exists it will be unique, and we say $N(p,q)$ is the hyperbolic $(p,q)$-Dehn filling of $\hat{N}$.  Thurston observed that when $\hat{N}$ has finite volume, the $(p,q)$-Dehn filling of $\hat{N}$ will exist for all but finitely many relatively prime pairs $(p,q)$ \cite{Thurston} (see also \cite{HK}).  Bonahon and Otal generalized this to geometrically finite manifolds \cite{BO} (see also \cite{Comar}).  Moreover, they show that if $(p_n, q_n)$ is an infinite sequence of distinct relatively prime pairs of positive integers, and $N(p_n, q_n)$ is the hyperbolic $(p_n, q_n)$-Dehn filling of $\hat{N}$ then $N(p_n,q_n)$ converges geometrically to $\hat{N}$.

%%%%
%%%% SECTION about punctured torus and review of Ken's notation
%%%%
\section{The Punctured Torus}   \label{bromberg section}

Recently, Ken Bromberg showed that when $\hat{T}$ is the punctured torus, $AH(\hat{T} \times I, \partial \hat{T} \times I)$ is not locally connected \cite{B}.  Specifically, he finds a geometrically finite representation $\rho$ with an additional parabolic and shows that near $\rho$, the deformation space is locally homeomorphic to a space $\mathcal{A}$ that is not locally connected.  Much of what follows depends on a detailed understanding of Bromberg's work. Accordingly, we plan to review the main results. 

Define $M_T = \hat{T} \times I$ and $P_T = \partial \hat{T} \times I$.  Let $\gamma$ be a nontrivial simple closed curve on $\hat{T}$ and define $P_T'$ to be the union of $P_T$ with an annular neighborhood of $\gamma$ on $\hat{T} \times \{1\}$.  The fundamental group of $M_T$ is a free group on two generators, which we will label $\pi_1(M_T) = \left< a,b \right>$.  Representations of $\pi_1(M_T)$ in $AH(M_T,P_T)$ are not only discrete and faithful, they also take the commutator $[a,b]$ to a parabolic element.  We can choose $\gamma$ in the homotopy class determined by $b$. 

Next we review the definition of the Maskit slice $\mathcal{M} \subset \mathbb{C}$.  If $\sigma \in R(M_T,P_T')$, then $\sigma : \pi_1(M_T) \to PSL(2, \mathbb{C})$ is conjugate into a one complex parameter family of representations determined by $z \in \mathbb{C}$.  Define $\sigma_z \in R(M_T, P_T)$ by   
\[
\sigma_z (a) = \begin{pmatrix}  iz & i \\ i & 0 \end{pmatrix} \; \text{and} \;  \sigma_z (b) = \begin{pmatrix}  1 & 2 \\ 0 & 1 \end{pmatrix}. 
\]
Of course not every $z \in \mathbb{C}$ yields a discrete, faithful representation $\sigma_z$.  We let 
\[
\mathcal{M} = \{z \in \mathbb{C} \; : \; \sigma_z \in AH(M_T, P_T') \}.
\]  
From the work of Keen, Series (section 2.5 of \cite{KS}) and Minsky \cite{Minsky} (summarized in section 4 of \cite{B}), $\mathcal{M}$ has two components in $\mathbb{C}$, one in the upper half plane, and its mirror image in the lower half plane (mirror in the sense of $z \mapsto -z$ and $z \mapsto \overline{z}$).  We will denote these by \[
\mathcal{M}^+ = \{ z \in \mathcal{M} \; : \; Im(z) > 0 \}  \; \text{and} \; \mathcal{M}^- = \{ z \in \mathcal{M} \; : \; Im(z) < 0 \}.
\]

Let $U$ be a regular neighborhood of $\gamma \times \{ \frac{1}{2} \}$ in $M_T$ and define $\hat{M} = M_T - U$.  Set $\hat{P} = P \cup \partial U$.  Express the fundamental group of $\hat{M}$ as $\pi_1(\hat{M}) = \left< a,b,c \: \vert \: [b,c] = 1\right>$.  Given $w \in \mathbb{C}$ we can extend a representation $\sigma_z \in MP(M_T, P_T')$ to a representation $\sigma_{z,w} \in R(\hat{M}, \hat{P})$ by defining
 \[
\sigma_{z,w} (a)  = \sigma_z (a), \; \sigma_{z,w} (b)  = \sigma_z (b), \; \sigma_{z,w} (c)  = \begin{pmatrix} 1 & w \\ 0 & 1 \end{pmatrix}.
\]

We may not always get a discrete, faithful representation, so for $z \in {\rm  int}(\mathcal{M}^+)$ we define the set
\[
\mathcal{A}_z = \{ w \in \mathbb{C} \; :\;  \sigma_{z,w} \in AH(\hat{M}, \hat{P}) \; \text{and}\; Im(w) > 0\},
\] and 
\[
\mathcal{A} = \{ (z,w) \in \mathbb{C} \times \hat{\mathbb{C}} \; : \; z \in {\rm  int}(\mathcal{M}^+),\; w \in \mathcal{A}_z \; \text{or}\; w = \infty\}.
\]

Then we have the following theorem of Bromberg:

\begin{thm} \label{bromberg1} If $z \in {\rm int}(\mathcal{M}^+)$, there is a neighborhood $V$ of $(z,\infty)$ in $\mathcal{A}$ and a continuous injective map $\Phi: V \to AH(M_T,P_T)$ such that $\Phi(V)$ contains a neighborhood of $\sigma_z$ in $AH(M_T,P_T)$. 
\end{thm}

We refer the reader to \cite{B} for a full discussion of how $\Phi$ is defined and what restrictions must be placed on $V$.  Using this local homeomorphism, Bromberg shows that $AH(M_T, P_T)$ is not locally connected by showing $\mathcal{A}$ is not locally connected.  We summarize what we will use from his work in the following theorem.

\begin{thm} \label{bromberg3}  There exists a point $z$ in the interior of $\mathcal{M}^+$ such that $\mathcal{A}$ is not locally connected at $(z,\infty)$.  Moreover, we can choose $w \in \mathbb{C}$ so that $(z,w + 2n) \in {\rm int}(\mathcal{A})$ for all integers $n$, $\Phi(z,w+2n)$ is defined and lies in $MP(M_T,P_T)$, and in any sufficiently small neighborhood $V$ of $(z, \infty)$ in $\mathcal{A}$, each point $(z, w + 2n)$ lies in a distinct component of $V$. 
\end{thm}

We now discuss how $\Phi$ is defined for the points $(z, w+2n)$ in the sequence above.  For points $(z,w) \in {\rm int}(\mathcal{A})$, we have $\sigma_{z,w} \in MP(\hat{M}, \hat{P})$ (see \cite{B}).  Let $\hat{N}_{z,w} = \mathbb{H}^3/ \sigma_{z,w}(\pi_1(\hat{M}))$ and let $\hat{f}_{z,w}: \hat{M} \to \hat{N}_{z,w}$ be the induced marking.  Let $\hat{f}_{z,w}(c)$ and $\hat{f}_{z,w}(b)$ mark the meridian and longitude of the rank $2$ cusp of $\hat{N}_{z,w}$ and define $N_{z,w}$ to be the hyperbolic $(1,0)$-filling of $\hat{N}_{z,w}$ (see Section \ref{filling notation} for this notation).  

Let \[
\phi_{z,w} : \hat{N}_{z,w} \to N_{z,w}
\]  denote the inclusion map.

Recall $z \in {\rm int}(\mathcal{M}^+)$ corresponds to a representation $\sigma_{z} \in MP(M_T,P_T')$, which in turn corresponds to a marked hyperbolic manifold $N_z = \mathbb{H}^3/\sigma_z (\pi_1(M_T))$, marked by a homotopy equivalence $f_z : M_T \to N_z$.  Since $\sigma_{z}(\pi_1(M_T))$ is a subgroup of $\sigma_{z,w}(\pi_1(\hat{M}))$, there is a covering map $$\pi_{z,w}: N_z \to \hat{N}_{z,w}.$$

The homotopy equivalence $f_{z,w} : M_T \to N_{z,w}$ defined by
\[
f_{z,w} = \phi_{z,w} \circ \pi_{z,w} \circ f_z.
\] defines a representation $(f_{z,w})_* : \pi_1(M_T) \to \pi_1(N_{z,w})$.  Note that as a hyperbolic manifold, we naturally identify $\pi_1(N_{z,w})$ with a subgroup of $PSL(2, \mathbb{C})$.   So we define $\Phi$ by 
\[
\Phi(z,w) = \begin{cases} (f_{z,w})_* \quad &\text{if} \; w \neq \infty \\ \sigma_{z} \quad &\text{if}\; w = \infty.  \end{cases}
\]

We make the following observation about the sequence of points $\Phi(z,w+2n)$ in Bromberg's construction.  

\begin{lem} \label{surgery lemma}  For a point $(z,w + 2n) \in {\rm int}(\mathcal{A})$ as in the sequence described above, the marked hyperbolic manifold $\Phi(z, w + 2n)$ is a $(1,n)$-Dehn filling of $\mathbb{H}^3/\sigma_{z,w}(\pi_1(\hat{M}))$. 
\end{lem}

\begin{proof}  A calculation shows $\sigma_{z,w+2n}(c) = \sigma_{z,w}(b^n c)$.  Thus $\pi_1(\hat{M})$ has the same image in $PSL(2, \mathbb{C})$ under all of the representations $\sigma_{z,w+2n}$, with the representations differing by automorphisms of $\pi_1(\hat{M})$ taking $a \mapsto a$, $b\mapsto b$ and $c \mapsto b^n c$.  This shows that a $(1,0)$-filling of $\hat{N}_{z,w +2n}$ is the same as a $(1,n)$-filling of $\hat{N}_{z,w}$.
\end{proof}

It now follows from Thurston's hyperbolic Dehn surgery theorem, generalized by Bonahon and Otal \cite{BO} to geometrically finite manifolds (see also Comar \cite{Comar}), that $N_{z,w+2n}$ converges geometrically to $\hat{N}_{z,w}$.  The following is a precise statement of this fact in terms of the geometric convergence of the corresponding Kleinian groups. 

\begin{lem} \label{geometric limit}   
There exists a convergent sequence of representations $\rho_n : \pi_1(M_T) \to PSL(2, \mathbb{C})$  such that $\rho_n$ is in the conjugacy class determined by $\Phi(z,w+2n) \in AH(M_T, P_T)$ and so that $\rho_n(\pi_1(M_T))$ converges geometrically to $\sigma_{z,w}(\pi_1(\hat{M}))$. 
\end{lem}

%%%%
%%%% SECTION  topological lemma
%%%%
\section{Detecting the Failure of Local Connectivity}   \label{topological section}

In Section \ref{main proof}, we will need a basic topological lemma to show that certain deformation spaces are not locally connected.

\begin{lem} \label{topological lemma}  Let $\Pi: X\to Y$ be a continuous map between topological spaces.  Let $U \subset Y$ be a neighborhood of $y_\infty \in Y$ such that there exists a sequence of points $y_n \to y_\infty$ in $U$ with each $y_n$ in a distinct component of $U$.   Suppose there exists a convergent sequence $x_n \to x_\infty$ in $X$ with $\Pi(x_n) = y_n$.  Then $X$ is not locally connected at $x_\infty$. 
\end{lem}

\begin{proof}
Let $U$ to be a neighborhood of $y_\infty$ as in the hypotheses.  Since $\Pi$ is continuous, $\Pi^{-1}(U)$ is an open set in $X$ that contains $x_\infty$.  If we suppose that $X$ is locally connected, then there must be an connected open neighborhood $W$ such that $x_\infty \in W \subset \Pi^{-1}(U)$.  By continuity, $\Pi(W) \subset U$ is connected.  However, since $x_n\to x_\infty$, $x_n \in W$ for all sufficiently large $n$, and since $\Pi(x_n) = y_n$, $\Pi(W)$ must be disconnected.  That is, $\Pi(W)$ is contained in $U$, but contains the points $y_n$ for all $n$ such that $x_n \in W$.  Since these are in distinct components of $U$, this contradicts that $\Pi(W)$ is connected.  
\end{proof}

%%%%
%%%% SECTION about first  main result
%%%%
\section{Constructing Manifolds whose Deformation Spaces are Not Locally Connected} \label{main proof}  

We are now ready to prove the main theorem.  Recall the notation $(M_T, P_T) = (\hat{T} \times I, \partial \hat{T} \times I )$. 

\textbf{Theorem \ref{main theorem}.} \textit{Let $M$ be a hyperbolizable 3-manifold containing a primitive essential annulus $A$, and suppose $(M_T, P_T)$ is pared homeomorphic to $(M', A)$, where $M'$ is the closure of one of the components of $M - A$.  If $P \subset \partial M$ is a paring locus that contains exactly one of the components of $\partial A$ and is otherwise disjoint from $M'$, then $AH(M,P)$ is not locally connected. }

\begin{proof}  By construction, there is a $\pi_1$-injective pared embedding $(M_T,P_T) \to (M,P)$ so any representation $\rho: \pi_1(M) \to PSL(2, \mathbb{C})$ restricts to a representation $\rho\vert_{\pi_1(M_T)}$ in $AH(M_T,P_T)$.  This restriction induces a continuous map
$$\Pi: AH(M,P) \to AH(M_T, P_T).$$  

Combining Theorems \ref{bromberg1} and \ref{bromberg3}, we can find a neighborhood $U$ of $\rho_\infty = \Phi(z,\infty)$ in $AH(M_T, P_T)$ and a sequence $\rho_n =  \Phi(z, w+2n)$ such that each $\rho_n $ lies in a distinct component of $U$. 
 
 We may conjugate $\rho_n$, $\rho_\infty$, and $\sigma_{z,w}$ such that the parabolic cyclic subgroup corresponding to the image of $[a,b]$ in $PSL(2, \mathbb{C})$ is generated by $\begin{pmatrix} 1 & 1 \\ 0 & 1 \end{pmatrix}$. We label the resulting Kleinian groups $G_n = \rho_n(\pi_1(N))$, $G_\infty = \rho_\infty(\pi_1(N))$, and $G = \sigma_{z,w}(\pi_1(\hat{N}))$.  Let $H$ denote the cyclic subgroup generated by $\begin{pmatrix} 1 & 1 \\ 0 & 1 \end{pmatrix}$.  Lemma \ref{geometric limit} shows that $G$ is the geometric limit of $G_n$.    
 
 To construct $\eta_n \in AH(M,P)$ with $\Pi(\eta_n) = \rho_n$ we need to apply Klein-Maskit combination using a uniform set of precisely invariant horoballs.  A set $B \subset \hat{\mathbb{C}}$ is precisely invariant under a subgroup $H$ in $G$ if (i) for all $h \in H$, $h(B) = B$ and (ii) for all $g \in G - H$, $g(B) \cap B = \emptyset$.  
 
Since $G$ is geometrically finite, there exists a pair of precisely invariant horoballs in $\Omega(G)$ for $H \subset G$ tangent to the point at $\infty$ in $\hat{\mathbb{C}}$, the fixed point of $\rho([a,b])$ (\textit{e.g.}, see p. 125 of \cite{Marden}).    We can take these to be 
  \[
 B_R^+ = \{ z \in \mathbb{C} \; : \; Im(z) > R \} \; \text{and} \; B_R^- = \{ z \in \mathbb{C} \; : \; Im(z) < -R \}
 \] for some $R$.  Without loss of generality, we can assume $R$ is large enough so that $\overline{B_R^\pm} \subset \Omega(G) - \{ \infty \}$.  Since $G_\infty$ is a subgroup of $G$ that contains $H$, the sets $B_{R}^\pm$ are precisely invariant under $H$ in $G_\infty$.  Next, we argue that these horoballs are precisely invariant for $H$ in $G_n$ for all but finitely many $n$.

 %%%%%
 %%%%% technical horoball lemma
 %%%%%
  \begin{lem} \label{invariant horoballs} The sets $B_R^+$ and $B_R^-$ (as above) are precisely invariant horoballs for $H$ in $G_n$, for all sufficiently large $n$.  
\end{lem}

\begin{proof} Our first claim is that the limit sets $\Lambda(G_n)$ are contained in the strip
\[
\Lambda(G_n) \subset  \{ z \in \mathbb{C} \; : \; |Im(z)| < R \} \cup \{ \infty \}
\]
for all sufficiently large $n$.

The groups $G_n$ are quasi-Fuchsian and the geometric limit $G$ is geometrically finite so the limit sets of $G_n$ converge to the limit set of $G$ in the Hausdorff topology on closed subsets of $\hat{\mathbb{C}}$ \cite{JM}, (see also Theorem 4.6.1 of \cite{Marden} or \cite{KT}).  Since $\Lambda(G_n) \to \Lambda(G)$ in the Hausdorff topology, the domains of discontinuity converge in the sense of Carath\'eodory. That is, if $K$ is a compact set in $\Omega(G)$ then $K \subset \Omega(G_n)$ for all sufficiently large $n$, and if $U$ is an open subset of infinitely many $\Omega(G_n)$ then $U \subset \Omega(G)$.  
Let 
\[
K = \{ z \; : \; 0 \leq Re(z) \leq 1, \; R \leq Im(z) \leq R + 1 \}.
\]  Clearly $K$ is compact and is contained in $\Omega(G)$.  Hence, $K \subset \Omega(G_n)$ for all sufficiently large $n$.  Since $H \subset G_n$ and $\Omega(G_n)$ is $G_n$-invariant, the entire strip
\[
W = \{ z \; :  \; R \leq Im(z) \leq R + 1 \} \subset \Omega(G_n)
\] is contained in $\Omega(G_n)$ for all sufficiently large $n$.

Since $G_n$ is quasi-Fuchsian, $\Lambda(G_n)$ is a Jordan curve, and by our choice of normalization, $\Lambda(G_n)$ goes through $\{\infty\}$ for all $n$.  Since $W \subset \Omega(G_n)$, the limit set of $G_n$ lies entirely above or below $W$.  That is, $\Lambda(G_n) - \{ \infty \}$ lies entirely above or below $W$.  The same argument can be applied to $-K$ to show that $-W \subset \Omega(G_n)$, and hence $\Lambda(G_n)  - \{ \infty \}$ is contained entirely in the upper half plane above $W$, entirely in the lower half plane below $-W$, or entirely in the desired strip $\{ z \in \mathbb{C} \; : \; |Im(z)| < R \}$.   Since $\Lambda(G)$ lies in this middle strip, and $\Lambda(G_n)$ converges to $\Lambda(G)$ in the Hausdorff topology, we must have $\Lambda(G_n)$ is entirely in this strip for all sufficiently large $n$.  

Now we address the precise invariance of $B_{R}^\pm$. If $h \in H$ then $h$ preserves horizontal lines so $B_R^\pm$ is invariant under $H$.  If $B_R^+$ is not precisely invariant under $H$, then there is an infinite sequence $g_{n_k} \in G_{n_k} - H$ such that $g_{n_k}(B_R^+) \cap B_R^+ \neq \emptyset$.  (The case $g_{n_k}(B_R^-) \cap B_R^- \neq \emptyset$ is similar.)  Without loss of generality, we can assume that for all $k$, $n_k$ is large enough so that $\Lambda(G_{n_k})$ is contained in the strip 
\[
\Lambda(G_{n_k}) \subset  \{ z \in \mathbb{C} \; : \; |Im(z)| < R \} \cup \{ \infty \}.
\]

If $g_{n_k}\left(\overline{B_R^+}\right) \subset \overline{B_R^+}$ then there is a fixed point of $g_{n_k}$ inside $\overline{B_R^+}$. This fixed point is in the limit set $\Lambda(G_{n_k})$, but $\Lambda(G_{n_k}) \cap \overline{B_R^+} = \{ \infty \}$ so $g_{n_k}$ must fix $\infty$.  However, $g_{n_k}$ cannot fix $\infty$ because it does not lie in the subgroup $H$, and $G_{n_k}$ is discrete. 

So we must have $g_{n_k}(B_R^+) \cap B_R^+ \neq \emptyset$, but $g_{n_k} \left(\overline{B_R^+}\right)$ is not contained in $\overline{B_R^+}$.  In this case, $g_{n_k}(\partial B_R^+) \cap \partial B_R^+$ contains some point in $\partial B_R^+  -  \{ \infty \}$. So there are points $x_{n_k} \in  \partial B_R^+  -  \{ \infty \}$ such that $g_{n_k}(x_{n_k}) \in \partial B_R^+  -  \{ \infty \}$.  Up to precomposition and postcomposition with elements of $H$, we can assume that $x_{n_k}$ and $g_{n_k}(x_{n_k})$ both lie in 
\[
L = \{ x + iR \; : \;  0 \leq x \leq 1 \} = \{ z \in \partial B_R^+ \; : \;  0 \leq  Re(z) \leq 1 \}. 
\] 
The set $L$ is compact so after passing to further subsequences if necessary $x_{n_k} \to x$ and $g_{n_k}(x_{n_k}) \to y$.  Then $g_{n_k}(x) \to y$.

Now by general convergence properties of Mobius transformations (see Theorem 2.1.1 of \cite{Marden}), we have one of two possibilities.  The first possibiilty is that, up to subsequence, $g_{n_k}$ converges to a Mobius transformation $g \in PSL(2, \mathbb{C})$.  By geometric convergence, $g \in G$.  This contradicts the fact that $B_R^+$ is precisely invariant for $H$ in $G$.  The alternative is the following:  suppose $p_{n_k}$ and $q_{n_k}$ are the repelling and attracting fixed points of $g_{n_k}$ (possibly equal to each other) with limits $p$ and $q$.    Then $g_{n_k}(z)$ converges to $q$ for $z \in \hat{\mathbb{C}} - \{p\}$.  But we have that $g_{n_k}(x) \to y$ so either $x = p$ or $y = q$.  Being limits of fixed points $p_{n_k}, q_{n_k}$, both $p,q$ lie in $\Lambda(G)$, but $R$ was chosen so that $\Lambda(G) \cap \overline{B_R^+} = \emptyset$.

It follows that $B_R^+$ is precisely invariant for $H$ in $G_n$ for all but finitely many $G_n$.  The same proof shows that $B_R^-$ is also precisely invariant for $H$ in $G_n$ for all but finitely many $n$. 
\end{proof}

%
%
%
 %%% lemma for construction of sequences
 %
 \begin{lem}  \label{sequence}
 There exists a convergent sequence of representations $\{\eta_n \}$ in $AH(M,P)$ such that $\Pi(\eta_n) = \rho_n$.
 \end{lem}
 
 \begin{proof}  
 By Lemma \ref{invariant horoballs}, we can find a fixed pair of horoballs $B_R^+$ and $B_R^-$ that are precisely invariant for $H$ in $G_n$ and for $H$ in $G_\infty$.  
 
 Recall that $M$ can be cut along the annulus $A$ resulting in two components, one of which is homeomorphic to $\hat{T} \times I$.  Call the other component $M_0$.  Note 
 \[
 \pi_1(M) = \pi_1(\hat{T}) *_{\pi_1(A)} \pi_1(M_0).
 \]   
 Let $P_0 = P \cap M_0$ and choose any $\eta_0 \in MP(M_0, P_{0})$.  Fix a representation in the conjugacy class determined by $\eta_0$ such that $\eta_0 (\pi_1(P_0)) = H$.  As $\eta_0$ is geometrically finite, there will be some $R'$ such that $B_{R'}^+$ and $B_{R'}^-$ is a pair of precisely invariant horoballs for $H$ in $\eta_0(\pi_1(M_0))$. Conjugate $\eta_0$ by $z \mapsto z + i (R + R')$ so the complement of $B_{R}^+$ is precisely invariant for $H$ in $\eta_0(\pi_1(M_0))$.  Define $\eta_n : \pi_1(M) \to PSL(2, \mathbb{C})$ by setting $\eta_n(g) = \rho_n(g)$ for $g \in \pi_1(\hat{T})$ and $\eta_n(g) = \eta_0(g)$ for $g \in \pi_1(M_0)$.  Similarly, define $\eta_\infty = \rho_\infty$ on $\pi_1(\hat{T})$ and $\eta_\infty = \eta_0$ on $\pi_1(M_0)$.   The first type of Klein-Maskit combination implies that $\eta_n, \eta_\infty \in AH(M,P)$ (\cite{AM}, \cite{Maskit original}, \cite{Maskit paper}).  Because $\rho_n \to \rho_\infty$ and the conjugacy representative of $\eta_0 \in MP(M_0, P_0)$ was fixed throughout the combination construction, $\eta_n \to \eta_\infty$.  It also follows that $\Pi(\eta_n) = \rho_n$.  
   \end{proof}
 
 We can now complete the proof of Theorem \ref{main theorem}. 
 The map $\Pi$ is continuous, and if we let $Y = AH(M_T,P_T)$, $y_n = \rho_n$, and $y_\infty = \rho_\infty$, then $Y$ contains a neighborhood $U$ and a sequence $y_n \to y_\infty$ satisfying the hypotheses of Lemma \ref{topological lemma}.  
  Lemma \ref{sequence} provides the sequences $x_n = \eta_n$ such that $\Pi(x_n) = y_n$. Hence $AH(M,P)$ is not locally connected at $\eta_\infty$. 
 
 \end{proof}

%% remarks at the end of first main section 

 \textit{Remark.}   Any geometrically finite representation in $AH(M_T, P_T)$ is in the image of $\Pi$. If $\rho \in AH(M_T, P_T)$, there is a conjugate of $\rho$ and an $R$ such that $B_{R}^\pm$ are precisely invariant under $H$ in $\rho(\pi_1(M_T))$.  So we can construct $\eta \in AH(M,P)$, as in Lemma \ref{sequence}, by applying Klein-Maskit combination to $\rho(\pi_1(M_T))$ and some group conjugate to $\eta_0 (\pi_1(M_0))$.  By construction, $\Pi(\eta) = \rho$. Note that the element conjugating $\eta_0$ so that Klein-Maskit combination can be applied may depend on $\rho$.

 Although the image of $\Pi$ contains all geometrically finite representations, the map is not surjective. An application of the covering theorem \cite{C} shows the image of $\Pi$ cannot contain a representation of the punctured torus whose limit set is all of $\hat{\mathbb{C}}$.

 %
 %
 % last section
 %  boundary connected-sum construction
 %
 \section{Another Family of Deformation Spaces that are Not Locally Connected} 
 Although we have considered manifolds $(M,P)$ which arise from gluing $\hat{T} \times I$ to another pared manifold along the annulus $\partial \hat{T} \times I$, our arguments can be adapted to show the following.

\textbf{Theorem \ref{second theorem}.} \textit{Suppose $M$ is obtained from a boundary connected sum of $\hat{T} \times I$ and a hyperbolizable manifold $(M',P')$ identifying a disk in $\hat{T} \times \{ 0 \}$ to a disk in $\partial M' - P'$.  Let $P = P' \cup (\partial \hat{T} \times I)$. Then $AH(M,P)$ is not locally connected. }

 Since the arguments are essentially the same as the previous section, we sketch the proof.  Again let $M_T = \hat{T} \times I$ and $P_T = \partial \hat{T} \times I$.  We still have a $\pi_1$-injective pared embedding of $(M_T, P_T)$ into $(M,P)$ and so we can define $\Pi : AH(M,P) \to AH(M_T, P_T)$ as before.  The points $\rho_n$ described in Lemma \ref{geometric limit} and their algebraic limit $\rho_\infty$ lie in the image of $\Pi$.   We express $\pi_1(M) = \pi_1(\hat{T}) * \pi_1(M')$, and choose any $\eta_0 \in MP(M',P')$.  We then construct $\eta_n : \pi_1(M) \to PSL(2, \mathbb{C})$ by choosing representations in the conjugacy classes determined by $\rho_n$ and $\eta_0$ (which we also denote $\rho_n$ and $\eta_0$) and setting $\eta_n(g) = \rho_n (g)$ for $g \in \pi_1(\hat{T})$ and $\eta_n(g) = \eta_0(g)$ for $g \in \pi_1(M')$.   
 
 We use the same normalization of $\rho_n$ and $\rho_\infty$ as in the beginning of the proof of Theorem \ref{main theorem}. Thus Lemma \ref{invariant horoballs} still applies, causing any disk $D \subset B_R^+$ of radius less than $1$ to be precisely invariant under the identity subgroup in $G_n$ (and $G_\infty$). One can find a fixed representation $\eta_0$ in the conjugacy class chosen above such that the complement of $D$ is precisely invariant under the identity in $\eta_0 (\pi_1(M'))$. It follows by Klein combination that $\eta_n, \eta_\infty \in AH(M,P)$. By construction $\eta_n \to \eta_\infty$ and $\Pi(\eta_n) = \rho_n$.  Thus, Lemma \ref{topological lemma} proves $AH(M,P)$ is not locally connected.

  \textit{Remark.} As noted in \cite{B}, Bromberg's work for $AH(\hat{T} \times I, \partial \hat{T} \times I)$ holds when $\hat{T}$ is replaced by the four-times punctured sphere $S_{0,4}$.  Thus one could obtain results analogous to Theorems \ref{main theorem}  and \ref{second theorem} by replacing $(M_T, P_T)$ with $(S_{0,4} \times I, \partial S_{0,4} \times I)$.

%%  font sizes from smallest to largest
%% \tiny
%% \scriptsize
%% \footnotesize
%% \small
%%  \normalsize
%% \large, \Large, \LARGE, \huge, \Huge

\footnotesize

%
%
%
%bibliography
%
%
%
%

\scriptsize

\textsc{Department of Mathematics \\ University of Michigan \\ 2074 East Hall \\ 530 Church St. \\ Ann Arbor, MI 48109 \\ magid@umich.edu}

\end{document}